\documentclass[12pt]{amsart}
\usepackage[all]{xy}
\usepackage{amsmath, amsfonts, amsthm , amssymb}
\usepackage{amscd, times, fullpage}
\newtheorem{proposition}{Proposition}[section]
\newtheorem{theorem}[proposition]{Theorem}

\theoremstyle{remark}
\newtheorem{defn}{Definition}[section]

\newtheorem{remark}[defn]{Remark}

\renewcommand{\L}{\ensuremath{\mathcal{L}}}
\newcommand{\Q}{\ensuremath{\mathbb{Q}}}

\newcommand{\Z}{\ensuremath{\mathbb{Z}}}

\newcommand{\lra}{\longrightarrow}
\newcommand{\Ker}{\mathop{\mathrm{Ker}}}

\newcommand{\Hom}{\ensuremath{\mathrm{Hom}}}
\newcommand{\rank}{\ensuremath{\mathrm{rank}}}

\begin{document}

%%% Title

\title[The based loop space on a flag manifold]{The integral Pontrjagin homology of the based loop space on a flag
manifold}
\author{Jelena Grbi\' c and Svjetlana Terzi\' c}
\address{School of Mathematics, University of Manchester, Manchester M13 9PL, United Kingdom}
\email{jelena.grbic@manchester.ac.uk}
\address{Faculty of Science, University of Montenegro, Podgorica, Montenegro}
\email{sterzic@cg.ac.yu}

\date{}

\begin{abstract}
The based loop space homology of a special family of homogeneous
spaces, flag manifolds of connected compact Lie groups is studied.
First, the rational homology of the based loop space on a complete
flag manifold is calculated together with its Pontrjagin structure.
Second, it is shown that the integral homology of the based loop
space on a flag manifold is torsion free. This results in a
description of the integral homology. In addition, the integral
Pontrjagin structure is determined.
\end{abstract}

\maketitle

\tableofcontents

\section{Introduction}

A \emph{complete flag manifold} of a compact connected Lie group $G$
is a homogeneous space $G/T$, where $T$ is a maximal torus in $G$.
In this paper we study the integral Pontrjagin homology of the based
loop space on a complete flag manifold $G/T$.

Compact homogeneous spaces, in particular, flag manifolds play a
significant role in many areas of physics and mathematics, such as
theory of characteristic classes of fibre bundles, representation
theory, string topology and quantum physics. Still there are only
few homogeneous spaces for which  the integral homology ring of
their based loop spaces is known. Some of them are classical simple
Lie groups, spheres, and complex projective spaces.

The motivation for our study comes from Borel's work~\cite{B} in
which he described the family of compact homogeneous spaces whose
cohomology ring is torsion free. In particular, homogeneous spaces
$G/U$ where $\rank\ G=\rank\ U$ stand out as homogeneous spaces
which behave nicely under application of algebraic topological
techniques. In this case Sullivan minimal model theory together with
the Milnor-Moore theorem can be employed to calculate the rational
homology ring of their based loop spaces. As one of the main results
of our paper (see Theorem~\ref{mainfree}) we prove that the homology
of the based loop space on a complete flag manifold is torsion free.

Furthermore, we explicitly calculate the integral Pontrjagin
homology ring of the loop spaces on the complete flag manifolds of
simple compact Lie groups $SU(n+1)$, $Sp(n)$, $SO(2n+1)$, $SO(2n)$, $G_{2}$,
$F_{4}$ and $E_6$ (see
Theorems~\ref{SUintegral},~\ref{Spintegral},~\ref{SO(2n+1)integral},
~\ref{SOintegral},~\ref{intG2/T},~\ref{intF4/T},~\ref{intE6/T}).

%Let $G$ be a compact connected Lie group, $U$ a closed connected
%subgroup and $G/U$ the corresponding homogeneous space. In this
%paper we study the homology of the based loop space on certain
%compact homogeneous spaces. There are only few homogeneous spaces
%for which  the integral homology ring of their based loop space is
%known. Some of them are: classical simple Lie groups, spheres and
%complex projective spaces.

%In these calculations we appeal to rational homotopy theory and to
%well known results on homology of the based loop space of simple
%compact Lie groups.

It is a classical result (see for example~\cite{BT}, or~\cite{BS})
that the homology of the $e$-connected component $\Omega _{0}G$ of
the loop space on $G$ is torsion free for any compact connected Lie
group $G$.
% and it is also known its homology algebra structure for
%simple compact Lie groups.
% and furthermore it is explicitly calculated for
Thus by the use of rational calculations, we show that there is a
split extension of algebras
\[
\xymatrix{ 1 \ar[r] & H_*(\Omega_{0}G; \Z)\ar[r] &
H_*(\Omega(G/T);\Z )\ar[r] &H_*(T;\Z )\ar[r] &1}
\]
%where $G$ is a simple compact Lie group and $T$ its maximal torus,
and
%consequently to
describe the integral Pontrjagin ring structure on $\Omega (G/T)$
for a simple compact Lie group $G$.

%of compact simple Lie groups where $\Omega G$ and $U$ are homologically
%torsion free, it may be possible to apply a spectral sequence
%argument to the fibration sequence $\ddr{\Omega
%G}{}{\Omega(G/U)}{}{U}$ in order to retrieve information on the
%existence of torsion in the homology of $\Omega (G/U)$. In the case
%when torsion does not exist, the rational homology determines the
%integral homology. On the other hand, when a classical Lie group has
%torsion (that is for $SO(m)$) another approach is needed. In
%particular, for $G=SO(n)$ and $U=T$ a maximal torus, we will consider
%certain fibrations involving homogeneous spaces $SO(2n+1)/T^n$ and $
%SO(2n)/T^n$, and show that $\Omega (SO(2n+1)/T^n)$ and $\Omega
%(SO(2n)/T^n)$ have torsion free homology.

%\begin{defn}
%A \emph{complete flag manifold} of a compact connected Lie group $G$
%is a homogeneous space $G/T$, where $T$ is a maximal torus in $G$.
%\end{defn}

%In this paper we calculate the integral homology of the loop space on
%a complete flag manifold of a classical compact Lie group and prove
%that it is torsion free.

%The motivation for our study comes from Borel's work~\cite{B} in
%which he described the family of compact homogeneous spaces whose
%cohomology ring is torsion free. As the main result of our paper
%(see Theorem~\ref{mainfree}) we prove  that the homology of the
%based loop space on a complete flag manifold is torsion free.

%Furthermore, we explicitly calculate the integral Pontrjagin
%homology ring of these loop spaces in
%Theorems~\ref{SUintegral},~\ref{Spintegral},~\ref{SO(2n+1)integral},~\ref{SOintegral}.

Throughout the paper, the loop space on a topological space will
mean a based loop space.

{\it Acknowledgements.} \ The authors would like to take this
opportunity to thank Professor Ralph Cohen for his helpful
suggestions and kind encouragement.

\section{Torsion in the homology of loop spaces}
We start by recalling some well known facts about the (co)homology
of classical simple compact Lie groups and their based loop spaces
(see for example~\cite{MT}). It is a classical result that for any
compact connected Lie group $G$ of rank $n$,
\[
H^{*}(G;\Q)\cong \wedge (z_1,\ldots ,z_n), \quad H_{*}(\Omega _{0}
G;\Q)\cong \Q [b_1,\ldots ,b_n]
\]
where $\deg (z_{i})=2k_{i}-1$ and $\deg (b_{i})= 2k_{i}-2$ for
$1\leq i\leq n$, and $k_{i}$ are  the exponents of the group $G$.
For simple compact Lie groups, these exponents are established.

For $G=SU(n+1)$ or $G=Sp(n)$, the integral homology of $G$ and
$\Omega G$ is torsion free and it is given by
\[
H^{*}(G;\Z)\cong \wedge (x_1,\ldots ,x_n), \quad H_{*}(\Omega
G;\Z)\cong \Z[y_1,\ldots ,y_n]
\]
where $\deg (x_{i})=2k_{i}-1$, and $\deg (y_{i})= 2k_{i}-2$ for
$1\leq i\leq n$. Under the rationalisation the integral generators
$x_1,\ldots, x_n$ and $y_1,\ldots, y_n$ are mapped onto the rational
generators $z_1,\ldots,z_n$ and $b_1,\ldots,b_n$, respectively.

For $G=SO(2n+1)$ or $G=SO(2n)$, the integral homology of $G$ and
$\Omega G$ has $2$-torsion.

Borel~\cite[Proposition~29.1]{B} proved that the homology of a flag
manifold $G/T$ is torsion free for the classical Lie groups $G$ and
for $G=G_2$ or $F_4$. Using Morse theory, it is proved in~\cite{BS1}
that this is true for any compact connected Lie group.

Our first result states that the complete flag manifold of a compact
connected Lie group behaves nicely with respect to the loop space
homology functor.
\begin{theorem}
\label{mainfree} The homology of the based loop space on the
complete flag manifold of a compact connected Lie group is torsion
free.
\end{theorem}
We will first show that to prove the theorem it is enough to
consider the case when $G$ is a simple, compact Lie group.
\begin{proposition}
\label{decompositionconnectedLie} The loop space on the flag
manifold of a compact, connected Lie group $G$ decomposes into a
product of the loop spaces on flag manifolds of simple, compact Lie
groups.
\end{proposition}
\begin{proof}
It is a classical result (see Onishchik~\cite{O}) that a compact
connected Lie group $G$ can be decomposed into a locally direct
product of connected simple normal subgroups. That is,
$G=G_{1}\cdots G_k$, where $G_i$ is a simple, connected Lie group or
a torus, $1\leq i\leq k$, such that
\[
\dim G_i\cap (G_1\cdots G_{i-1}\cdot G_{i+1}\cdots G_k) = 0.
\]
Let $\widetilde{G}$ be $G_1\times \ldots \times G_k$ and $p\colon
\widetilde{G}\to G$ defined by $p(g_1,\ldots
,g_k)=g_1\cdots g_k$. Since $\Ker p =
\cup_{i=1}^{n}G_i\cap (G_1\cdots \hat{G_i}\cdots G_k)$, we obtain
that $\Ker p$ is discrete or in other words $p\colon
\widetilde{G}\to G$ is a covering. Thus $\Ker p$ is contained in the
center $Z(\widetilde{G})$ of $\widetilde{G}$. Let $T=T_1\times
\ldots \times T_k$ be a maximal torus in $\widetilde{G}$, where
$T_i$ is a maximal torus in $G_i$ for $1\leq i\leq k$. Then $\Ker
p\subset T$ and therefore
\[
G_1/T_1\times \ldots\times G_k/T_k = \widetilde{G}/T =
(\widetilde{G}/\Ker p)/T = G/T.
\]
Hence
\[
\Omega (G/T)\simeq \Omega (G_1/T_1)\times \ldots\times \Omega
(G_k/T_k).
\]
\end{proof}

{\it Proof of Theorem~\ref{mainfree}.} Let $G$ be a compact
connected Lie group and $T$ its maximal torus. We have that the
complete flag manifold $G/T$ for any compact connected Lie group $G$
is homeomorphic to the complete flag manifold $\widetilde{G}/T$ of
its universal cover $\widetilde{G}$. Therefore, we may assume $G$ to
be simply connected. For $G$ simply connected, it is classical
result (see for example~\cite{S}) that  $\Omega (G/T)$ has the same
homotopy type as $\Omega (G)\times T$. To verify this notice that
related to the principal fibration  
as topological spaces. For $G$ a simple, compact, simply connected Lie group, it is a classical
result that the integral homology of $\Omega G$  is torsion free. Now using splitting
\eqref{simplyconnectedsplitting}, we conclude that the homology of
$\Omega (G/T)$ is torsion free in this case. The statement of the theorem now
follows readily from Proposition~\ref{decompositionconnectedLie}.\qed
%When $G$ is not simply connected, the inclusion $i\colon T\lra G$ is
%not null homotopic so the above splitting does not exist and
%moreover $\Omega G$ is not torsion free. On the other hand, we have that
%the complete flag manifold $G/T$ for any compact connected Lie
%group $G$ is homeomorphic to the complete flag manifold of its
%universal cover $\widetilde{G}$. Now as there is a homotopy decomposition
%$\Omega (G/T)\cong\Omega(\widetilde{G}/T)\simeq
%\Omega\widetilde{G}\times T$, the homology of $\Omega
%(G/T)$ is torsion free.

%This generalises to any connected group.

%{\it Proof of Theorem~\ref{mainfree}.} The statement of the theorem
%follows readily from Proposition~\ref{decompositionconnectedLie}  and
%the preceding observation that the homology of the loop space on
%the flag manifold of a simple compact Lie group is torsion free.\qed

\section{Rational homology}
In this section we calculate the rational homology ring of the loop
space on a flag manifold by looking separately at each simple Lie
group.

To calculate the rational homology of the based loop space on a
complete flag manifold of a classical simple Lie group we will apply
Sullivan minimal model theory. Let us start by recalling the key
constructions and setting the notation related to the Sullivan
minimal model and rational homology of loop spaces which we are
going to use in the subsequent sections.

\subsection{Rational homology of loop spaces}
Let $M$ be a simply connected topological space with the rational
homology of finite type. Let $\mu=(\Lambda V, d)$ be a  Sullivan
minimal model for $M$. Then $d\colon V\lra \Lambda^{\geq 2}V$ can be
decomposed as $d=d_1+d_2+\cdots$, where $d_i\colon V\lra
\Lambda^{\geq i+1}V$. In particular, $d_1$ is called the
\emph{quadratic part} of the differential $d$.

The homotopy Lie algebra $\L$ of $\mu$ is defined in the following
way. Define a graded vector space $L$ by requiring that
\[
sL=\Hom(V, \Q)
\]
where as usual the suspension $sL$ is defined by $(sL)_i=(L)_{i-1}$.
We can define a pairing $\langle \, ; \, \rangle\colon V\times sL
\lra \Q$ by $\langle v;sx\rangle =(-1)^{\deg v}sx(v)$ and extend it
to $(k+1)$-linear maps
\[
\Lambda^kV\times sL\times\cdots\times sL\lra \Q
\]
by letting
\[
\langle v_1\wedge\cdots\wedge v_k; sx_k,\ldots
,sx_1\rangle=\sum_{\sigma\in S_k}\epsilon_\sigma\langle
v_{\sigma(1)}; sx_1\rangle\cdots\langle v_{\sigma(k)};sx_k\rangle
\]
where $S_k$ is the symmetric group on $k$ letters and
$v_{\sigma(1)}\wedge\cdots\wedge v_{\sigma(k)}=\epsilon_\sigma
v_1\wedge\cdots\wedge v_k$. It is important to notice that $L$
inherits a Lie bracket $[\, ,
 \, ]\colon L\times L\lra L$ from $d_1$ uniquely determined by
 \begin{equation}
\label{bracket}
 \langle v; s[x,y]\rangle= (-1)^{\deg y+1}\langle
 d_1v;sx,sy\rangle \quad \text{ for } x,y\in L, v\in V.
 \end{equation}
Denote by $\L$ the Lie algebra $(L, [\, ,\, ])$.

Recall that the graded Lie algebra $L_M=(\pi_*(\Omega M)\otimes \Q;
[\, , \,])$ is called the {\it rational homotopy Lie algebra of
$M$}. The commutator $[\, , \,]$ is given by the Samelson product.
There is an isomorphism between the rational homotopy Lie algebra
$L_M$ and the homotopy Lie algebra $\L$ of $\mu$. Using the
theorem in the Appendix of Milnor and Moore~\cite{MM}, it follows
that
\[
H_*(\Omega M; \Q)\cong U\L
\]
where $U\L$ is the universal enveloping algebra for $\L$. Further
on,
\[
U\L \cong T(L)/\langle xy-(-1)^{\deg x\deg y}yx-[x,y]\rangle.
\]
For a more detailed account of this construction see for
example~\cite{FHT}, Chapters  $12$ and $16$.

As the notion of formality will be important for our calculation we
recall it here.
\begin{defn}
A commutative cochain algebra $(A,d)$ satisfying $H^0(A)=\Q$ is {\it
formal} if it is weakly equivalent to the cochain algebra
$(H(A),0)$.
\end{defn}
Thus $(A,d)$ and a path connected topological space $X$ are formal
if and only if their minimal Sullivan models can be computed
directly from their cohomology algebras.
\begin{remark}\label{borelideal}
There are  some known cases of topological spaces for which a
minimal model can be explicitly computed and  formality proved. Some
of them, that are important for us in this work, are the spaces that
have so called ``good cohomology" in terminology of~\cite{BG}.
Namely, topological space $X$ is said to have good cohomology if
\[
H^*(X; \Q )\cong \Q [u_1,\ldots , u_n]/\langle P_1,\ldots , P_k
\rangle
\]
where the polynomials $P_1,\ldots ,P_k$ form the regular sequence
in $\Q [u_1, \ldots ,u_n]$, or in other words, the ideal
$\langle P_1,\ldots , P_k\rangle$ is a Borel ideal in $\Q [u_1,
\ldots ,u_n]$. In this case Bousfield and Gugenheim~\cite{BG} proved
that the minimal model of $X$ is given by
\[
\mu (X) = \Q [u_1, \ldots , u_n]\otimes \wedge (v_1, \ldots, v_k)
\]
where $\deg (v_i)=\deg (P_i)-1$ for $1\leq i\leq k$, and the
differential $d$ is given by
\[
d(u_i) = 0, \quad d(v_j) = P_j.
\]
\end{remark}

\subsection{The loop space on a complete flag manifold}
In this section we calculate the rational homology of
the loop space on the complete flag manifold of a simple Lie
group.

%Borel~\cite{B} proved that the homology of a homogeneous space $G/U$,
%where $\rank G=\rank U$ is torsion free. Since the rank of $SU(n+1)$
%is $n$, the homology  $H^*(SU(n+1)/T^n;\Z)$ is torsion free.
Recall from Borel~\cite[Section $26$]{B} that the rational (as well
as integral) cohomology of $SU(n+1)/T^n$ is the polynomial algebra
on $n+1$ variables of degree 2 quotient out by the ideal generated
by the symmetric functions in these variables
\[
H^*(SU(n+1)/T^n;\Q)\cong \Q [u_1,\ldots ,u_{n+1}]/\langle
S^{+}(u_1,\ldots , u_{n+1}) \rangle .
\]
It is important to note that the ideal $\langle S^{+}(u_1,\ldots ,
u_{n+1}) \rangle$ is a Borel ideal. As a consequence, by
Remark~\ref{borelideal}, $SU(n+1)/T^n$ is formal. Thus the minimal
model for $SU(n+1)/T^n$ is the  minimal model for the commutative
differential graded algebra $(H^*(M;\Q), d=0)$ and it is given by
$\mu=(\Lambda V, d)$, where
\[
V=(u_1,\ldots ,u_{n},v_1,\ldots ,v_n) \] and $\deg (u_{k})=2$,
$\deg(v_k)=2k+1$ for $1\leq k\leq n$.

The differential $d$ is defined by
\begin{equation}\label{differential}
d(u_{k})=0, \quad d(v_{k})=\sum _{i=1}^{n}u_{i}^{k+1} +
(-1)^{k+1}\left ( \sum _{i=1}^{n}u_{i}\right )^{k+1}.
\end{equation}
It is easy to see that the quasi isomorphism $f\colon\mu =(\Lambda V,
d)\lra (H^*(M;\Q), d=0)$ is given by the following rule
\[
u_{i}\mapsto u_{i}, \quad v_{i}\mapsto 0  \text{ for } 1\leq i\leq n
.
\]
\begin{theorem}
\label{flagQ} The rational homology ring of the loop space on the flag
manifold $SU(n+1)/T^n$ is
\begin{equation}
\label{QSU} H_*(\Omega ( SU(n+1)/T^n); \Q) \cong
\end{equation}
\[
\indent \left( T(a_1,\ldots ,a_n)/\left\langle
a_{k}^2=a_pa_q+a_qa_p\, |\, 1\leq k,p,q\leq n, p\neq q \right\rangle\right)
\otimes \Q[b_2,\ldots ,b_n]
\]
where the generators $a_{i}$ are of degree $1$ for $1\leq i\leq n$,
and the generators $b_{k}$ are of degree $2k$ for $2\leq k\leq n$.
\end{theorem}

\begin{proof}
The underlying vector space of the homotopy Lie algebra $\L$ of
$\mu$ is given by
\[
L=(a_1,\ldots ,a_n,b_1,\ldots ,b_n)\] where $\deg(a_k)=1,\
\deg(b_k)=2k$ for $1\leq k\leq n$.

In order to define Lie brackets we need the quadratic part $d_1$ of
the differential in the minimal model. In this case, using the
differential $d$ defined in~\eqref{differential}, the quadratic part
$d_1$ is given by
\[
d_1(u_l)=0  \text{ for } 1\leq l\leq n,   \ d_1(v_1)=2\sum
_{i=1}^{n} u_i^2 + 2\sum _{i<j}u_iu_j,\ d_{1}(v_k)=0 \ \text{ for }
k\neq 1.
\]

For dimensional reasons, we have
\[
[a_k, b_l]=[b_k, b_l]=0 \text{ for } 1\leq k,l\leq n.
\]
By the defining property of the Lie bracket stated
in~\eqref{bracket}, we have

\[
\langle v_1; s[a_k,a_k]\rangle=\left\langle 2\sum u_i^2 + 2\sum
u_iu_j; sa_k, sa_k\right\rangle =2\langle u_k^2; sa_k,sa_k\rangle =4
\text{ and }
\]
\[
\langle v_1; s[a_k,a_l]\rangle= 2\langle u_ku_l; sa_k,sa_l\rangle = 2
\text{ for } k\neq l \
\]
resulting in the commutators
\[
[a_k,a_l]=2b_1 \text{ for } k\neq l, \text{ and } [a_k,a_k]=4b_1.
\]
Therefore in the tensor algebra $T(a_1,\ldots , a_n,b_1,\ldots
,b_n)$, the Lie brackets above induce the following relations
\[\begin{array}{ll}
a_ka_l+a_la_k =2b_1 & \text{ for }  1\leq k,l\leq n, k\neq l, \\
a_k^2=2b_1 & \text{ for }  1\leq k\leq n,\\
a_kb_l=b_la_k &  \text{ for }  1\leq k,l\leq n, \\
b_kb_l=b_lb_k &  \text{ for }  1\leq k,l\leq n.
\end{array}
\]
Thus
\begin{equation}
 U\L  \cong \left( T(a_1,\ldots ,a_n)/\left\langle a_k^2 =a_pa_q+a_qa_p
  \right\rangle \right)\otimes \Q[b_2,\ldots ,b_n].
\end{equation}
This proves the theorem.
\end{proof}

The rational cohomology rings for the flag manifolds
$SO(2n+1)/T^n\cong Spin(2n+1)/T^n$, $SO(2n)/T^n\cong Spin(2n)/T^n$, and
$Sp(n)/T^n$ (see for example Borel~\cite[Section $26$]{B}) are given
by
\[
H^{*}(SO(2n+1)/T^n;\Q)\cong H^{*}(Sp(n)/T^n;\Q)\cong\Q [u_1,\ldots
,u_n]/\langle S^{+}(u_1^2,\ldots ,u_n^2) \rangle ,
\]
\[
H^{*}(SO(2n)/T^n;\Q) \cong \Q [u_1,\ldots ,u_n]/\langle
S^{+}(u_1^2,\ldots ,u_n^2),u_1\cdots u_n\rangle
\]
where $u_i$ is of degree 2 for $1\leq i\leq n$.

By Remark~\ref{borelideal}, all the above mentioned complete flag manifolds
are formal and therefore their minimal Sullivan model is the minimal
model for their cohomology algebra with the trivial differential.

Proceeding in the same way as in the previous theorem, we obtain the
following results.
\begin{theorem}
\label{rationalOSO(2n+1)} The rational homology ring of the loop space on
$SO(2n+1)/T^n$ and $Sp(n)/T^n$ is given by
\begin{equation}
\label{QSp} H_*(\Omega (SO(2n+1)/T^n); \Q)\cong H_*(\Omega
(Sp(n)/T^n); \Q) \cong
\end{equation}
\[
\left( T(a_1,\ldots ,a_n)/\left\langle \begin{array}{l}a_1^2=\ldots =a_n^2,\\
a_ka_l=-a_la_k \text{ for } k\neq l\end{array}\right\rangle
\right)\otimes \Q[b_2,\ldots ,b_n]
\]
where the generators $a_i$ are of degree $1$ for $1\leq i\leq n$,
and the generators $b_k$ are of degree $4k-2$ for $2\leq k\leq n$.
\end{theorem}
\begin{proof}
We give just an outline of the proof as it is similar to the proof
of Theorem~\ref{flagQ}. The minimal model for $SO(2n+1)/T^n$ is
given by $\mu=(\Lambda V, d)$, where
\[
V=(u_1,\ldots ,u_{n},v_1,\ldots ,v_n),
\]
and $\deg (u_{k})=2,$ $\deg (v_k)=4k-1$  for $1\leq k\leq n$.

The differential $d$ is given by
\begin{equation}\label{differentialSO(2n+1)}
d(u_{k})=0, \quad d(v_{k})=\sum _{i=1}^{n}u_{i}^{2k}\ \text{ for }
1\leq k\leq n.
\end{equation}
Therefore the underlying vector space of the homotopy Lie algebra
$\L$ of $\mu$ is
\[
L=(a_1,\ldots ,a_n,b_1,\ldots ,b_n)
\]
where $\deg (a_{k})=1,$ $\deg (b_k)=4k-2$  for $1\leq k\leq n$, and
the quadratic part $d_1$ of the differential $d$ is given by
\[
d_1(u_l)=0 \text{ for } 1\leq l\leq n, \ d_1(v_1)=\sum _{i=1}^{n}
u_i^2,\ d_{1}(v_k)=0 \ \text{ for }  k\geq 2.
\]
The induced Lie brackets on $L$ are equal to
\[
\begin{array}{ll}
 \text{[}a_k, b_l]=[b_k, b_l]=0 & \text{ for } 1\leq k,l\leq n,\\
\text{[}a_k, a_k]=2b_1 &\text{  for } 1\leq k\leq n,\\
\text{[}a_k, a_l]=0 &\text{ for } k\neq l.
\end{array}
\]
This implies the following relations in $U\L$:
\[
\begin{array}{ll}
a_{k}^{2}=b_1 &\text{ for } 1\leq k\leq n,\\
a_ka_l+a_la_k=0 & \text{ for } k\neq l,\\
a_kb_l=b_la_k & \text{ for }  1\leq k,l\leq n, \\
b_kb_l=b_lb_k & \text{ for }  1\leq k,l\leq n.
\end{array}
\]
The theorem follows now at once knowing that $H_*(\Omega
(SO(2n+1)/T^n); \Q)\cong U\L$.
\end{proof}
\begin{theorem}
\label{rationalOSO(2n)} The rational homology ring of the loop space on
$SO(2n)/T^n$ for $n>2$ is given by
\[
H_*(\Omega (SO(2n)/T^n); \Q)\cong
\]
\[
\left( T(a_1,\ldots ,a_n)/\left\langle
\begin{array}{l}
a_1^2=\ldots =a_n^2,\\
a_ka_l=-a_la_k \text{ for } k\neq l
\end{array}\right\rangle
\right)\otimes \Q[b_2,\ldots ,b_{n-1},b_n]
\]
where the generators $a_i$ are of degree $1$ for $1\leq i\leq n$, the
generators $b_k$ are of degree $4k-2$ for $2\leq k \leq n-1$, and the
generator $b_n$ is of degree $2n-2$.
\end{theorem}
\begin{proof}
To be reader friendly we outline a proof. The minimal model for
$SO(2n)/T^n$ is given by $\mu=(\Lambda V, d)$, where
\[
V=(u_1,\ldots ,u_{n},v_1,\ldots ,v_{n-1},v_n),
\]
and $\deg (u_{k})=2,$ $\deg (v_k)=4k-1$  for $1\leq k\leq n-1$ and
$\deg (v_n)=2n-1$.

The differential $d$ is given by
\begin{equation}\label{differentialSO(2n)}
d(u_{k})=0, \quad d(v_{k})=\sum _{i=1}^{n}u_{i}^{2k} \text{ and }
d(v_n)=u_1\cdots u_n.
\end{equation}
Hence the underlying vector space of the homotopy Lie algebra $\L$
of $\mu$ is
\[ L=(a_1,\ldots ,a_n,b_1,\ldots ,b_{n-1},b_n)
\]
where $\deg (a_{k})=1,$ $\deg (b_k)=4k-2$  for $1\leq k\leq n-1$,
$\deg (b_n)=2n-2$, and the quadratic part $d_1$ of the differential
$d$ is given by
\[
d_1(u_l)=0 \text{ for } 1\leq l\leq n, \ d_1(v_1)=\sum _{i=1}^{n}
u_i^2,\ d_{1}(v_k)=0 \ \text{ for } 2\leq k\leq n.
\]
The induced Lie brackets on $L$ are equal to
\[
\begin{array}{ll}
 \text{[}a_k, b_l]=[b_k, b_l]=0 &\text{ for } 1\leq k,l\leq n,\\
\text{[}a_k, a_k]=2b_1 &\text{ for } 1\leq k\leq n, \\
\text{[}a_k, a_l]=0 &\text{ for } k\neq l,
\end{array}
\]
and thus in $U\L$:
\[
\begin{array}{ll}
a_{k}^{2}=b_1 &\text{ for } 1\leq k\leq n,\\
a_ka_l+a_la_k=0 & \text{ for } k\neq l,\\
a_kb_l=b_la_k & \text{ for }  1\leq k,l\leq n, \\
b_kb_l=b_lb_k & \text{ for }  1\leq k,l\leq n.
\end{array}
\]
Since $H_*(\Omega(SO(2n)/T^n); \Q)\cong U\L$, we have proved the
theorem.
\end{proof}

In the theorems that follow we compute the rational homology rings
of the based loop space on the complete flag manifolds of the
exceptional Lie groups $G_2$, $F_4$ and $E_6$. We refer to~\cite{C}
and~\cite{M} for the Weyl group invariant polynomials which we use
for the descriptions of the rational cohomology rings of the
complete flag manifolds of these groups.  We want also to emphasize
that the rational, as well as the integral, cohomology rings of the
flag manifolds $G_2/T^2$, $F_4/T^4$ and $E_6/T^6$ are thoroughly
discussed in~\cite{TW}.

\begin{theorem}
The rational homology ring of the loop space on $G_2/T^2$ is given
by
\[
H_{*}(\Omega (G_2/T^2); \Q )\cong \left(T(a_1,a_2)/\langle
a_1a_2+a_2a_1=a_1^2=a_2^2 \rangle\right)\otimes \Q [b_5]
\]
where $\deg b_5 = 10$, and $\deg a_1 = \deg a_2 =1$.
\end{theorem}
\begin{proof}
Recall that
\[
H^{*}(G_2/T^2;\Q )\cong \Q [u_1,u_2,u_3]/\langle P_1,P_2,P_6\rangle
\]
where $P_1=\sum\limits_{i=1}^{3}u_1$, $P_2=\sum\limits_{i=1}^{3}u_i^2$,
$P_6=\sum\limits_{i=1}^{3}u_i^6$ and
$\deg u_1=\deg u_2=\deg u_3=2$. Therefore the minimal model is
$\Lambda V=\Lambda(u_1,u_2,v_1,v_5)$ where $\deg v_1=3$,  $\deg v_5
=11$, and the differential $d$ is given by $d(u_1)=d(u_2)=0,
d(v_1)=2(u_1^2+u_2^2+u_1u_2), d(v_5)=u_1^6+u_2^6+(u_1+u_2)^6$. Thus
\[
d_1(u_1)=d_1(u_2)=d_1(v_5)=0, \; d_1(v_1)=2(u_1^2+u_2^2+u_1u_2).
\]
In the homotopy Lie algebra $L=(a_1,a_2,b_1,b_5)$ the induced
commutator relations are given by
\[
[a_i,b_j]=0 \text{ for } i=1,2, \ j=1,5,\; [b_i,b_j]=0 \text{ for }
i,j = 1,5, [a_1,a_2]=2b_1, \; [a_1,a_1]=[a_2,a_2]=4b_1.
\]
Hence the following relations in $U\L$ hold:
\[
\begin{array}{l}
 a_ib_j-b_ja_i=0  \text{ for }  i=1,2, j=1,5,\\
 b_1b_5=b_5b_1,\\
a_1a_2+a_2a_1=2b_1,\\
 a_1^2=a_2^2=2b_1.
\end{array}
\]
\end{proof}
\begin{theorem}
The rational homology ring of the loop space on $F_4/T^4$ is given
by
\[
H_{*}(\Omega(F_4/T^4);\Q )\cong
\]
\[
\left(T(a_1,a_2,a_3,a_4)/\left\langle
\begin{array}{l}
a_1^2=\ldots=a_4^2,\\ a_ia_j=-a_ja_i\ for\
 i\neq j\end{array}\right\rangle\right)\otimes \Q [b_5,b_7,b_{11}]
\]
where $\deg a_i=1$ for $1\leq i\leq 4$, $\deg b_5=10$, $\deg b_7=14$,
and $\deg b_{11}=22$.
\end{theorem}
\begin{proof}
The rational cohomology algebra of $F_4/T^4$ is
\[
H^{*}(F_4/T^4;\Q ) \cong \Q [u_1,u_2,u_3,u_4]/\langle
P_2,P_6,P_8,P_{12} \rangle
\]
where $\deg u_i =2$ for $1\leq i\leq 4$, and
\[
P_{k}=u_1^k+u_2^k+u_3^k+u_4^k + \frac{1}{2^{k+1}}(\pm u_1\pm u_2\pm u_3\pm
u_4)^k
\]
for $k=2,6,8,12$. For degree reasons, the only relevant generator
for determining $d_1$ is \\$P_2 = 3(u_1^2+u_2^2+u_3^2+u_4^2)$.
Therefore we have
\[
V = (u_1,u_2,u_3,u_4,v_2,v_6,v_8,v_{12}), \mbox{ where } \deg
v_k=2k-1
\]
and the quadratic part of $d$ is given by
\[
d_1(u_i)=0 \text{ for } 1\leq i\leq 4,\; d_1(v_j)=0 \text{ for }
j=6,8,12,
\]
\[
\text{ and } d_1(v_2)=3(u_1^2+u_2^2+u_3^2+u_4^2).
\]
This determines the homotopy Lie algebra
\[
L = (a_1,a_2,a_3,a_4,b_1,b_5,b_7,b_{11})
\]
where $\deg a_i=1,\; \deg b_1=2,\; \deg b_5=10,\; \deg b_7=14$, and
$\deg b_{11}=22$ with the Lie brackets given by
\[
\begin{array}{ll}
 \text{[}a_i,b_j]=[b_l,b_j]=0 & \text{ for } 1\leq j\leq 4 \text{ and } j,l=1,5,7,11,\\
 \text{[}a_i,a_j]=0 & \text{ for } i\neq j, \\
 \text{[}a_i,a_i]=6b_1 &\text{ for } 1\leq i\leq 4.
\end{array}
\]
This implies that in $U\L$ for every possible $i$ and $j$, $a_i$ and
$b_j$ commute as well as $b_i$ and $b_j$ does. Also the additional
relations in $U\L$ hold:
\[
a_1^2=a_2^2=a_3^2=a_4^2=3b_1, \text{ and } a_ia_j+a_ja_i=0 \text{
for } i\neq j.
\]
The statement of the theorem now follows directly.
\end{proof}
\begin{theorem}
The rational homology ring of the loop space on $E_6/T^6$ is given
by
\[
H_{*}(\Omega(E_6/T^6);\Q)\cong
\]
\[\indent \left( T(a_1,\ldots ,a_5,a)/\left\langle
\begin{array}{l}
a^2=a_{k}^2=a_pa_q+a_qa_p\ for\ 1\leq k, p,q\leq
5,p\neq q \\
aa_i=-a_ia\ for\ 1\leq i\leq 5\end{array}\right\rangle\right) \otimes \Q[b_4,b_5,b_7,b_8,b_{11}],
\]
where $\deg a_i=1$ for $1\leq i\leq 5$, $\deg a=1$,  and $\deg b_j=2j$ for
$j=4,5,7,8,11$.
\end{theorem}
\begin{proof}
The rational cohomology of $E_6/T^6$ is
\[
H^{*}(E_6/T^6;\Q )\cong
\Q [u_1,u_2,u_3,u_4,u_5,u_6, u]/\langle
P_1,P_2,P_5,P_6,P_8,P_9,P_{12} \rangle
\]
where $\deg u_i=2$ for $1\leq i\leq 6$, $\deg u=2$, and
\[
P_{k}=\sum_{i=1}^{6}(u_i\pm u)^k +\sum_{1\leq i<j\leq 6}(-1)^{k}(u_i+u_j)^k
\]
for $k=2,5,6,8,9,12$, and $P_{1}=\sum\limits_{i=1}^{6}u_i$. It follows that
$V=(u_1,u_2,u_3,u_4,u_5,u,v_2,v_5,v_6,v_8,v_9,v_{12})$ and $d_1$
is determined only by
\[
P_2=12(u_1^2+\ldots +u_5^2+u^2 +\sum_{1\leq i< j\leq 5}u_iu_j).
\]
In a similar fashion as before we obtain that
$L=(a_1,a_2,a_3,a_4,a_5,a,b_1,b_4,b_5,b_7,b_8,b_{11})$, where $\deg
a_i=\deg a=1 $ for $1\leq i\leq 5$, and $\deg b_j=2j$ for
$j=1,4,5,7,8,11$. The commutators are
\[
[a_i,a]=[a_i,b_j]=[a,b_j]=[b_l,b_j]=0 \text{ for } 1\leq i\leq 5 \text{ and } j,l=1,4,5,7,8,11,
\]
\[
[a_i,a_j]=12b_1 \text{ for } 1\leq i,j\leq 5, i\neq j, \text{ and }
[a_i,a_i]=[a,a]=24b_1  \text{ for } 1\leq i\leq 5.
\]
The last three commutator relations imply the following relations in
$U\L$:
\[
a^2=a_i^2=a_ka_l+a_la_k=12b_1 \text{ for } 1\leq i,l,k\leq 5, k\neq l.
\]
This directly implies the statement of the theorem.
\end{proof}

\section{Integral Pontrjagin homology}
%Borel~\cite[Proposition~29.1]{B} proved that the homology of a complete flag
%manifold $G/T$ is torsion free for classical $G$ and for $G=G_2$ or $F_4$. Using Morse theory,
%it is proved in~\cite{BS1} that this is true for any compact connected Lie group.
In this section we study the integral Pontrjagin ring structure of
$\Omega (G/T)$, where $G$ is a simple Lie group. We make use of the
rational homology calculations for $\Omega (G/T)$ from the previous
section and the results from~\cite{BT},~\cite{M} and~\cite{W} on
integral homology of the identity component $\Omega _{0}G$ of the
loop space on $G$. Recall that $H_{*}(\Omega _{0}G;\Q)$ is
primitively generated for a compact connected Lie group $G$.

\subsection{The integral homology of $\Omega (SU(n+1)/T^n)$}
\begin{theorem}
\label{SUintegral} The integral Pontrjagin homology ring of the loop
space on $SU(n+1)/T^n$ is
\[
H_*\left(\Omega (SU(n+1)/T^n); \Z\right)\cong
\]
\[
\left( T(x_1,\ldots ,x_n)\otimes \Z[y_1,\ldots
,y_n]\right)/\left\langle x_k^2=x_px_q+x_qx_p=2y_1 \text{ for }
1\leq k,p,q\leq n, p\neq q\right\rangle
\]
where the generators $x_1,\ldots ,x_n$ are of degree $1$, and the
generators $y_i$ are of degree $2i$ for $1\leq i\leq n$.
\end{theorem}
\begin{proof}
It is well known that if $G$ is a simply connected Lie group, then
$\pi_{2}(G/T)\cong \Z ^{\dim T}$ and $\pi _{3}(G/T)\cong \Z$. Let
\[
W\colon \pi _{2}(G/T)\otimes \pi _{2}(G/T)\to \pi _{3}(G/T)
\]
denote the pairing given by the Whitehead product. In what follows,
we identify $H_{1}(T,\Z )$ with $\pi _{2}(G/T)$ and $H_{2}(\Omega
G,\Z )$ with $\pi _{3}(G/T)$ via natural homomorphisms. Thus since
there is no torsion in homology, and using the rational homology
result~\eqref{QSU}, we obtain that there is a split extension of
algebras
\[
\xymatrix{ 1 \ar[r] & H_*(\Omega SU(n+1); \Z)\ar[r] &
H_*(\Omega(SU(n+1)/T^n);\Z )\ar[r] &H_*(T^n;\Z )\ar[r] &1}
\]
with the extension given by $[\alpha ,\beta]=W(\alpha ,\beta)\in H_{2}(\Omega SU(n+1);\Z )$, where
$\alpha ,\beta \in H_{1}(T^n;\Z )$.

We explain the extension of the algebra in more detail. Notice that
there is a monomorphism of two split extensions of algebras
\[
\xymatrix{ 1 \ar[r]\ar[d] & H_*(\Omega SU(n+1); \Z)\ar[r]\ar[d] &
H_*(\Omega(SU(n+1)/T^n);\Z )\ar[r]\ar[d] &H_*(T^n;\Z )\ar[r]\ar[d] &
1\ar[d]\\
1 \ar[r]& H_*(\Omega SU(n+1); \Q)\ar[r] & H_*(\Omega(SU(n+1)/T^n);\Q
)\ar[r] &H_*(T^n; \Q)\ar[r]&1.}
\]
Denote by $\bar{c}_2,\ldots ,\bar{c}_{n+1}$ the universal
transgressive generators in $H^{*}(SU(n+1);\Z )$ which map to the
symmetric polynomials $c_2=\sum\limits_{1\leq i<j\leq
n+1}x_ix_j,\ldots ,c_{n+1}=x_1\cdots x_nx_{n+1}$ generating
$H^{*}(BSU(n+1);\Z)$. The elements  $x_1,\ldots ,x_n, x_{n+1}$ are
the integral generators of $H_{*}(T^n; \Z)$ and
$\sum\limits_{i=1}^{n+1}x_i=0$. Now let $y_1,\ldots ,y_n$ be the
integral generators of $H_{*}(\Omega SU(n+1);\Z )$ obtained by the
transgression of the elements from $H_{*}(SU(n+1);\Z )$ which are
the Poincare duals of $\bar{c}_2,\ldots ,\bar{c}_{n+1}$. Further,
the subspace of primitive elements in $H_{*}(\Omega SU(n+1);\Z )$ is
spanned by the elements $\sigma _1,\ldots ,\sigma _n$ which can be
expressed in terms of $y_1,\ldots, y_n$ using the Newton formula
\begin{equation}\label{Newton}
\sigma
_{k}=\sum_{i=1}^{k-1}(-1)^{i-1}\sigma_{k-i}y_{i}+(-1)^{k-1}ky_{k},
\;\; 1\leq k\leq n.
\end{equation}
The integral elements $\sigma _1,\ldots ,\sigma _n$ rationalise to
the elements $b_1,\ldots ,b_n\in H_{*}(\Omega SU(n+1);\Q )$. The
generators $a_1,\ldots,a_n$ in $H_*(T^n;\Q )$ are the rationalised
images of the integral generators $x_1,\ldots, x_n$ in $H_{*}(T^n;\Z
)$. To decide the integral extension, we consider the rational
Pontrjagin ring structure~\eqref{QSU} of $\Omega (SU(n+1)/T^n)$.
Looking at the  above commutative diagram of the algebra extensions,
we conclude that the integral elements
\[\begin{array}{ll}
x_kx_l+x_lx_k -2\sigma_1 & \text{ for }  1\leq k,l\leq n, k\neq l, \\
x_k^2-2\sigma_1 & \text{ for }  1\leq k\leq n,\\
x_k\sigma_l-\sigma_lx_k &  \text{ for }  2\leq k,l\leq n, \\
\sigma_k\sigma_l-\sigma_l\sigma_k &  \text{ for }  2\leq k,l\leq n
\end{array}
\]
from $H_{*}(\Omega(SU(n+1)/T^n);\Z )$ map to zero in
$H_*(\Omega(SU(n+1)/T^n);\Q)$. As the map between the algebra
extensions is a monomorphism, we conclude that these integral
elements are zero. Using that there is no torsion in homology and
Newton formula~\eqref{Newton}, we have
\[
\begin{array}{ll}
x_kx_l+x_lx_k=2y_1 & \text{ for } 1\leq k,l\leq n, k\neq l, \\
x_k^2=2y_1 & \text{ for }  1\leq k\leq n,\\
x_ky_l-y_lx_k=0&  \text{ for }  2\leq k,l\leq n, \\
y_ky_l-y_ly_k=0 &  \text{ for }  2\leq k,l\leq n
\end{array}
\]
which completely describes the integral Pontrjagin ring of
$\Omega(SU(n+1)/T^n)$ and finishes the proof.
\end{proof}

\subsection{The integral homology of $\Omega (Sp(n)/T^n)$}

\begin{theorem}
\label{Spintegral} The integral Pontrjagin homology ring of the
based loop space on $Sp(n)/T^n$ is
\[
H_*(\Omega (Sp(n)/T^n); \Z)\cong \]
\[
\left( T(x_1,\ldots ,x_n)/\left\langle
\begin{array}{l}
x_1^2=\ldots =x_n^2,\\
x_kx_l=-x_lx_k \text{ for } k\neq l
\end{array}
\right\rangle\right)\otimes \Z [y_2,\ldots ,y_n]
\]
where the generators $x_1,\ldots ,x_n$ are of degree $1$, and the
generators $y_i$ are of degree $4i-2$ for $2\leq i\leq n$.
\end{theorem}
\begin{proof}
The proof is analogous to the proof of Theorem~\ref{SUintegral}.
Denote by $\bar{c}_1,\ldots ,\bar{c}_{n}$ the universal
transgressive generators in $H^{*}(Sp(n);\Z )$ which map to the
generators $c_1=\sum\limits_{i=1}^{n}x_i^2,\ldots ,c_{n}=x_1^2\cdots
x_n^2$ of $H^{*}(BSp(n);\Z )$. Let $y_1,\ldots, y_n$ be the
generators in $H_{*}(\Omega Sp(n);\Z)$ obtained by the transgression
of the elements in $H_{*}(Sp(n);\Z )$ which are the Poincare duals
of $\bar{c}_1,\ldots ,\bar{c}_{n}$. Recall from~\cite{BT} that the
subspace of the primitive elements in $H_{*}(\Omega Sp(n);\Z )$ is
spanned by the elements $\sigma _1,\ldots ,\sigma _n$ given by
\begin{equation}\label{NS}
\sigma _{k}=\sum\limits_{i=1}^{k-1}(-1)^{i-1}\sigma_{k-i}
y_{i}+(-1)^{k-1}ky_{k},\;\; 1\leq k\leq n.
\end{equation}
The integral elements $\sigma_1,\ldots, \sigma _n$ rationalise to
the generators $b_1,\ldots ,b_n$ of $H_{*}(\Omega Sp(n);\Q )$ given
in~\eqref{QSp}.  The generators $a_1,\ldots,a_n$ of $H_{*}(T^n;\Q )$
are the rationalised images of the integral generators $x_1,\ldots,
x_n$ in $H_*(T^n;\Z)$ . Therefore we conclude that in
$H_{*}(\Omega(Sp(n)/T^n);\Z )$ the following  integral elements are
zero:
\[\begin{array}{ll}
x_kx_l+x_lx_k -\sigma_1 & \text{ for }  1\leq k,l\leq n, k\neq l \\
x_k^2-\sigma_1 & \text{ for }  1\leq k\leq n\\
x_k\sigma_l-\sigma_lx_k &  \text{ for }  2\leq k,l\leq n \\
\sigma_k\sigma_l-\sigma_l\sigma_k &  \text{ for }  2\leq k,l\leq n
\end{array}.
\]
Since there is no torsion in homology, going back to Newton
formula~\eqref{NS}, we obtain the same relations between $y_1,\ldots
,y_n$ and $x_1,\ldots, x_n$ which determine the integral Pontrjagin
ring structure on $\Omega (Sp(n)/T^n)$.
\end{proof}

\subsection{The integral homology of $\Omega (SO(2n)/T^n)$ and
$\Omega (SO(2n+1)/T^n)$ }

As mentioned before, $SO(m)$ is  not simply connected and the
cohomology of $SO(m)$ and the homology of $\Omega SO(m)$ are not
torsion free, namely, they have $2$-torsion. Nevertheless, since
$SO(m)/T\cong Spin(m)/T$, where $T$ is a maximal torus, the rational
homology calculations enable us to prove the following.
\begin{theorem}
\label{SO(2n+1)integral} The integral Pontrjagin homology ring of
the based loop space on $SO(2n+1)/T^n$ is given by
\[
H_*(\Omega (SO(2n+1)/T^n); \Z)\cong \left( T(x_1,\ldots ,x_n)
\otimes\Z [y_1,\ldots,y_{n-1},2y_n,\ldots 2y_{2n-1}]\right)/I
\]
where $I$ is  generated by
\[
\begin{array}{l}
x_1^2-y_1,\;\; x_i^2-x_{i+1}^2\; \text{ for }\;  1\leq i\leq n-1 \\
x_kx_l+x_lx_k\; \text{ for }\; k\neq l\\
y_{i}^2 - 2y_{i-1}y_{i+1}+\ldots \pm 2y_{2i}\; \text{ for }\; 1\leq
i\leq n-1
\end{array}
\]
where $\deg x_i =1$ for $1\leq i\leq n$,  $\deg y_{i}=2i$ for $1\leq
i\leq 2n-1$, $\deg 2y_{i}=2i$ for $n\leq i\leq 2n-1$, and $y_0=1$.
\end{theorem}

\begin{remark}\label{BTSO(2n+1)}
Before proving  Theorem~\ref{SO(2n+1)integral}, let us  recall the
ring structure of $H_{*}(\Omega _{0}SO(2n+1); \Z)$. It is proved
in~\cite{BT} that the algebra  $H_{*}(\Omega _{0}SO(2n+1); \Z )$ is
generated by the classes $y_1,\ldots ,y_{n-1}, 2y_n,\ldots
,2y_{2n-1}$ which satisfy the relations
\[
y_{i}^2 - 2y_{i-1}y_{i+1}+2y_{i-2}y_{i+2}-\ldots \pm 2y_{2i}=0
\text{ for } 1\leq i\leq n-1
\]
where $\deg y_i=2i$ for $1\leq i\leq n-1$, $\deg 2y_{i}=2i$ for
$n\leq i\leq 2n-1$, and $y_0=1$. For  $[\frac{n+1}{2}]\leq i\leq
n-1$, these relations express $2y_{2i}$ in terms of $y_1,\ldots,
y_{n-1}, 2y_n,\ldots ,2y_{2i-1}$  and thus eliminate  $2y_{2i}$ as
generators. For $1\leq i\leq [\frac{n+1}{2}]-1$, the relations above
imply new relations on the generators $y_{2i}$, that is,
$2y_{2i}=\pm (y_{i}^2-2y_{i-1}y_{i+1}+\ldots \pm 2y_1y_{2i-1})$.
This implies that the elements $y_{2i}$ for $1\leq i\leq
[\frac{n+1}{2}]-1$ are generators only in the homology of
$\Omega_{0}SO(2n+1)$ with coefficients where $2$ is not invertible.
Consider the rational elements $p_k$ defined by the recursion formula
\begin{equation}
  \label{recursion}
p_{k}-p_{k-1}y_1+\ldots \pm ky_k=0 \text{ for } 1\leq k\leq 2n-1
\text{ where } p_{0}=1.
\end{equation}
The relations in $H_{*}(\Omega _{0}SO(2n+1);\Z )$ imply that only
$p_1,p_3,\ldots, p_{2n-1}$ are non zero. According to~\cite{BT} the
elements $p_1,p_3,\ldots, p_{2[\frac{n}{2}]-1},
2p_{2[\frac{n}{2}]+1}, \ldots, 2p_{2n-1}$ span the subspace of
primitive elements in $H_{*}(\Omega _{0} SO(2n+1);\Z)$. These
elements are obtained by transgressing the elements in
$H_{*}(SO(2n+1);\Z )$ which are the Poincare duals of the universal
transgressive generators $\bar{\sigma}_1,\ldots ,\bar{\sigma}_n$ in
$H^{*}(SO(2n+1);\Z )$. The generators $\bar{\sigma}_1,\ldots
,\bar{\sigma}_n$ map to the symmetric polynomials $\sigma
_{i}(x_1^2,\ldots ,x_n^2)$ for $1\leq i\leq n$ generating the free
part in $H^{*}(BSO(2n+1);\Z)$. In this way we see that
$p_1,p_3,\ldots, 2p_{2[\frac{n}{2}]+1},\ldots ,2p_{2n-1}$
rationalise to the rational generators $b_{i}$ in $H_{*}(\Omega _{0}
SO(2n+1);\Q )$ (see Theorem~\ref{rationalOSO(2n+1)}).
\end{remark}

\begin{remark}
If we denote the generators of $H_{*}(\Omega _{0}SO(2n+1);\Z )$ by
$y_1,\ldots,y_{n-1},y_n,\ldots, y_{2n-1}$, then the relations are
slightly more complicated and they are given by
\[
y_{i}^2+
2\sum_{k=1}^{min\{i,n-1-i\}}(-1)^ky_{i-k}y_{i+k}+
\sum_{k=n-i}^{i}(-1)^ky_{i-k}y_{i+k}=0
\]
where $1\leq i\leq n-1$.
\end{remark}
\begin{proof}
Recall that $SO(2n+1)/T^n\cong Spin(2n+1)/T^n$ implying that $\Omega
(SO(2n+1)/T^n)\cong \Omega (Spin (2n+1)/T^n)$. It is known that
$\Omega Spin(2n+1)\cong\Omega _{0}SO(2n+1)$, see for
example~\cite{MT}. Consider the morphism of two extensions of
algebras
\[
\xymatrix{H_*(\Omega _{0}SO(2n+1); \Z)\ar[r]\ar[d] &
H_*(\Omega(SO(2n+1)/T^n);\Z )\ar[r]\ar[d] &H_*(T^n;\Z )\ar[d]\\
H_*(\Omega_0SO(2n+1); \Q)\ar[r] & H_*(\Omega(SO(2n+1)/T^n);\Q
)\ar[r] &H_*(T^n; \Q).}
\]
By Remark~\ref{BTSO(2n+1)}, we have that all the generators
$b_1,\ldots, b_n$ of $H_*(\Omega _{0}SO(2n+1);\Q )$ are in the
rationalisation of the integral elements $p_1,p_3,\ldots,
p_{2[\frac{n}{2}]-1},2p_{2[\frac{n}{2}]+1},\ldots ,2p_{2n-1}$ of
$H_*(\Omega _{0}SO(2n+1); \Z)$. Since the map between two algebra
extensions is a monomorphism, we conclude that in
$H_{*}(\Omega(SO(2n+1)/T^n);\Z )$ the following relations hold
\[
\begin{array}{ll}
x_kx_l+x_lx_k=p_1 & \text{ for }  1\leq k,l\leq n, k\neq l, \\
x_k^2=p_1 & \text{ for }  1\leq k\leq n,\\
x_kp_{2l-1}=p_{2l-1}x_k &  \text{ for }  2\leq k,l\leq n, \\
p_{2k-1}p_{2l-1}=p_{2l-1}p_{2k-1} &  \text{ for }  2\leq k,l\leq n
\end{array}
\]
as these elements map to zero in $H_{*}(\Omega(SO(2n+1)/T^n);\Q )$.
Note that $p_1=y_1$, which gives that $y_1=x_1^2$ in $H_{*}(\Omega
(SO(2n+1)/T^n);\Z)$.

The fact that differs this case from the case of $SU(n+1)$ or
$Sp(n)$ is that these integral elements $p_1,\ldots ,2p_{2n-1}$ that
map onto rational generators, do not produce all the generators in
$H_{*}(\Omega _{0}SO(2n+1);\Z)$. Nevertheless, since there is no
torsion in homology, we can also deduce from the rational homology
calculations that there is a split extension of algebras
\[
\xymatrix{ 1 \ar[r] & H_*(\Omega_{0}SO(2n+1); \Z)\ar[r] &
H_*(\Omega(SO(2n+1)/T^n);\Z )\ar[r] &H_*(T;\Z )\ar[r] &1.}
\]
We have that $y_{2i-1}$ survive as the generators
in $H_{*}(\Omega(SO(2n+1)/T^n);\Z)$ for $2\leq i\leq n$ using
the relations coming from $H_{*}(\Omega _{0}SO(2n+1);\Z)$ and the fact
that the integral elements $p_3,\ldots ,2p_{2n-1}$ rationalise to the
generators $b_2,\ldots ,b_n$ in
$H_{*}(\Omega(SO(2n+1)/T^n);\Q)$.
Therefore, in order to verify the above splitting we need
to show that the generators $y_{2i}$ for $1\leq i\leq [\frac{n+1}{2}]$ in
$H_{*}(\Omega _{0}SO(2n+1); \Z )$ survive as  generators in $H_{*}(\Omega(SO(2n+1)/T^n); \Z)$.
We prove this by
induction on $i$. If $y_2$ is not a generator in
$H_{*}(\Omega(SO(2n+1)/T^n); \Z)$, then it can be expressed as
\[
y_{2}=\alpha x_{1}^{4} + \sum\limits_{i=2}^{n} \beta
_{i}x_{1}^{3}x_{i}+\sum\limits_{2\leq i< j\leq n}\gamma
_{ij}x_{1}^{2}x_{i}x_{j}+\sum\limits_{1\leq i<j<k<l\leq n}\delta
_{ijkl}x_{i}x_{j}x_{k}x_{l},
\]
where $\alpha ,\beta_i ,\gamma _{ij}, \delta_{ijkl}$ are  integers.
On the other hand, in $H_{*}(\Omega _{0}SO(2n+1); \Z)$ we have that
$2y_2=y_1^2$ which translates to $2y_2=x_1^4$ in
$H_{*}(\Omega(SO(2n+1)/T); \Z)$. This implies that $\beta
_{i}=\gamma_{ij}=\delta_{ijkl}=0$, and $2\alpha =1$, which is
impossible since $\alpha$ is an integer.
In the same way, assuming
that $y_{2i}$ for $1\leq i\leq k<[\frac{n+1}{2}]$ are generators in
$H_{*}(\Omega(SO(2n+1)/T^n); \Z)$, we prove that $y_{2(k+1)}$ is a
generator as well. If it were not, we would have
\[
y_{2(k+1)}=\alpha y_{k+1}^2 +P(x_1,\ldots,x_n,y_2,\ldots,y_{2k+1})
\]
where $\alpha\in \Z$ and $P$ is a polynomial with integer
coefficients which does not contain $y_{k+1}^2$. On the other hand,
in the relation $2y_{2(k+1)}=\pm (y_{k+1}^2-2y_{k}y_{k+2}-\cdots \pm
2y_{1}y_{2k+1})$ in $H_{*}(\Omega _{0} SO(2n+1); \Z)$, when translating
to $H_{*}(\Omega(SO(2n+1)/T^n); \Z)$ we have by the inductive
hypothesis that $y_{k+1}^2$ can not be eliminated. This implies that
the coefficient $\alpha$ satisfies $2\alpha =\pm 1$ which is
impossible.

%Recall that these elements are  generators in $H_{*}(\Omega_{0}(SO(2n+1), \F )$ for any field $\F$ whose  %characteristic is not $2$. When $F = \Z _{2}$, it is known that $H_{*}(\Omega_{0}(SO(2n+1)), \Z _{2}) =
%\Z _{2}[y_1,\ldots ,y_n]$ in degrees less than $2n$. Therefore in degrees greater than $2n-1$ we may have in %$H_{*}(\Omega _{0}(SO(2n+1), \Z )$ new generators $y$ and $2y$ are going to be polynomials in $y_1,\ldots ,y_n$ with %integral coefficients.

We are left with a verification of the commutator relations in
$H_{*}(\Omega (SO(2n+1)/T^n);\Z )$. Since $2y_2=x_{1}^4$, we have
$2y_2x_i=x_{1}^4x_i=x_ix_{1}^4=2x_iy_2$, that is, $y_2x_i=x_iy_2$.
Now by induction on $k$, we prove that $y_{k}x_j=x_jy_{k}$ for an
arbitrary $y_{k}$. For $k$  odd,  relation~\eqref{recursion}
together with the inductive hypothesis gives that  $x_i$ for  $1\leq
i\leq n$ commutes with $y_k$. Let $k$ be even. Since $\deg y_{i}$ is
even for any $i$, each monomial in the polynomial
$P(x_1,\ldots,x_n,y_2,\ldots,y_{k-1})=2y_k$  contains even number of
generators $x_1,\ldots,x_n$. Using now the inductive hypothesis, we
have that every $x_i$ commutes with $P$ and thus with $y_k$.
\end{proof}

\begin{theorem}
\label{SOintegral} The integral Pontrjagin homology ring of the
based loop space on $SO(2n)/T^n$ is given by
\[
H_*(\Omega (SO(2n)/T^n); \Z)\cong
\]
\[
\left( T(x_1,\ldots ,x_n)\otimes \Z [y_1,\ldots, y_{n-2},y_{n-1}+z,
y_{n-1}-z,2y_n,\ldots ,2y_{2(n-1)}]\right)/I
\]
where $I$ is generated by
\[
\begin{array}{l}x_1^2-y_1,x_i^2-x_{i+1}^2\; \text{ for }\; 1\leq i
\leq n-1\\
x_kx_l+x_lx_k\; \text{ for }\; k\neq l\\
y_i^2-2y_{i-1}y_{i+1}+2y_{i-2}y_{i+2}-\ldots \pm 2y_{2i}\;
 \text{ for }\; 1\leq i\leq n-2\\
(y_{n-1}+z)(y_{n-1}-z)-2y_{n-1}y_{n+1}+\ldots \pm 2y_{2(n-1)}
\end{array}
\]
where $\deg x_i =1$ for $1\leq i\leq n$, $\deg y_i=2i$ for $1\leq
i\leq n-2$, $\deg(y_{n-1}+z)=\deg(y_{n-1}-z)=2(n-1)$, $\deg
2y_{i}=2i$ for $n\leq i\leq 2(n-1)$ and $y_0=1$.

%the
%generators $y_k$ are of degree $4k-2$ for $2\leq k \leq n-1$ and the
%generator $y_n$ is of degree $2n-2$.
\end{theorem}
\begin{remark}
Recall from~\cite{BT} that the algebra $H_{*}(\Omega _{0}SO(2n);\Z)$
is generated by the elements
$y_1,\ldots,y_{n-2}, y_{n-1}+z,y_{n-1}-z, 2y_{n},\ldots ,2y_{2(n-1)}$ which satisfy the relations
\[
y_i^2-2y_{i-1}y_{i+1}+2y_{i-2}y_{i+2}-\ldots \pm 2y_{2i} =0
\text{ for } 1\leq i\leq n-2,
\]
\[
(y_{n-1}+z)(y_{n-1}-z)-2y_{n-1}y_{n+1}+\ldots \pm 2y_{2(n-1)}=0.
\]
As in previous case, these relations eliminate $2y_{2i}$ as
generators for $[\frac{n+1}{2}]\leq i\leq n-2$, while for $1\leq
i\leq [\frac{n+1}{2}]-1$, they induce new relations on $y_{2i}$
implying that $y_{2i}$ are generators only in the homology of
$\Omega_{0}SO(2n)$ with coefficients where $2$ is not invertible.
The subspace of primitive elements in $H_{*}(\Omega _{0}SO(2n);\Z )$
is spanned by the elements $p_1,p_3,\ldots p_{n-2},2z,2p_{n},\ldots
,2p_{2(n-1)-1}$ for $n$ odd and by the elements
$p_1,p_3,\ldots,p_{n-1},2z,2p_{n+1},\ldots ,2p_{2(n-2)+1}$ for $n$
even. These primitive generators are obtained by transgressing the
elements in $H_{*}(SO(2n);\Z )$ which are the Poincare duals of the
universal transgressive generators $\bar{\sigma}_1,\ldots
,\bar{\sigma}_{n-1},\bar{\lambda}$  in $H^{*}(SO(2n);\Z )$. The generators
$\bar{\sigma}_1,\ldots ,\bar{\sigma}_{n-1},\bar{\lambda}$ map to the polynomials
$\sigma _{i}(x_1^2,\ldots ,x_n^2)$ for $1\leq i\leq n-1$ and
$\lambda =x_1\cdots x_n$ which generate the free part in
$H^{*}(BSO(2n);\Z)$.
\end{remark}

\begin{proof}
The proof is analogous to the proof of
Theorem~\ref{SO(2n+1)integral}.
\end{proof}
\subsection{The integral homology of $\Omega (G_2/T^2)$}
\begin{theorem}\label{intG2/T}
The integral Pontrjagin homology ring of $\Omega (G_2/T^2)$ is given
by
\[
H_{*}(\Omega (G_2/T^2); \Z )\cong \left(T(x_1,x_2)\otimes
\Z[y_1,y_2,y_5]\right)/\langle x_1x_2+x_2x_1=x_1^2=x_2^2=2y_1,
2y_2=x_1^4 \rangle,
\]
where $\deg x_1=\deg x_2=1$, $\deg y_2=4$, and $\deg y_5 =10$.
\end{theorem}
\begin{remark}\label{intG2}
The integral homology algebra of $\Omega G_2$ has the following
form~\cite{BT}:
\[
H_{*}(\Omega G_2; \Z ) = \Z [y_1,y_2,y_5]/\langle 2y_2-y_1^2\rangle,
\]
with $\deg y_1=2$, $\deg y_2=4$, and $\deg y_5=10$.
\end{remark}
\begin{proof}
Consider a morphism  of two  extensions of algebras
\[
\xymatrix{ H_*(\Omega G_2; \Z)\ar[r]\ar[d] &
H_*(\Omega(G_2/T^2);\Z )\ar[r]\ar[d] &H_*(T^2;\Z )\ar[d]\\
 H_*(\Omega G_2; \Q)\ar[r] & H_*(\Omega(G_2/T^2);\Q
)\ar[r] &H_*(T^2; \Q).}
\]
In $H_*(\Omega G_2;\Q )$ the generators $b_1$ and $b_5$ are the
rationalisations of the integral elements $y_1$ and $y_5$ in
$H_*(\Omega G_2; \Z)$. It follows that the relations between $x_1,x_2$ and
between $y_1,y_5$ in $H_{*}(\Omega (G_2/T^2);\Z)$ are lifted from the relations on their
rationalisations. We further show that there is a split extensions
of algebras
\[
\xymatrix{ 1 \ar[r] & H_*(\Omega G_2; \Z)\ar[r] &
H_*(\Omega(G_2/T^2);\Z )\ar[r] &H_*(T^2;\Z )\ar[r] &1.}
\]
To deduce the splitting above, we use that there is no torsion in
the corresponding homologies. We first need to show that the
generator $y_2\in H_{*}(\Omega G_2; \Z )$ survives as a generator in
$H_{*}(\Omega(G_2/T^2); \Z)$.  If it were not, we would have that
$y_2=\alpha x_1^4 + \beta x_1^3x_2$, and that $2y_2=x_1^4$ using the
relations in $H_{*}(\Omega G_2; \Z )$. This would imply that $2\alpha = 1$
which is impossible since $\alpha$ is an integer. Since $2y_2=x_1^4$
and there is no torsion in homology, using already established
relations, we get that $y_2$ commutes with other generators in
$H_{*}(\Omega(G_2/T^2); \Z)$.
\end{proof}
\subsection{The integral homology of $\Omega (F_4/T^4)$}
\begin{theorem}\label{intF4/T}
The integral Pontrjagin homology ring of $\Omega (F_4/T^4)$ is given
by
\[
H_{*}(\Omega (F_4/T^4); \Z )\cong
\]
\[
\left(T(x_1,x_2,x_3,x_4)\otimes \Z
[y_1,y_2,y_3,y_5,y_7,y_{11}]\right)/I
\]
where $I=\langle x_i^2=3y_1,\; 1\leq i\leq 4,\; x_ix_j=x_jx_i,\;
i\neq j,\; 2y_2=x_1^4,\; 3y_3=x_1^2y_2\rangle$, where $\deg x_i=1$
for $1\leq i\leq 4$, and $\deg y_i=2i$ for $i=2,3,5,7,11$.
\end{theorem}
\begin{remark}
The integral homology algebra $H_{*}(\Omega F_4;\Z )$ is computed
in~\cite{W} and it is given by
\[
H_{*}(\Omega F_4;\Z ) = \Z [y_{1},y_{2},y_{3},y_{5},y_{7},y_{11}]/
\langle y_{1}^2-2y_{2},y_{1}y_{2}-3y_{3}\rangle.
\]
\end{remark}
\begin{proof}
As in the previous cases, we first prove that there is a split
extension of algebras
 \[
\xymatrix{ 1 \ar[r] & H_*(\Omega F_4; \Z)\ar[r] &
H_*(\Omega(F_4/T^4);\Z )\ar[r] &H_*(T^4;\Z )\ar[r] &1.}
\]
Since there is no torsion in homology, the rational homology
calculations for $\Omega (F_4/T^4)$ gives that it is enough to prove
that $y_2$ and $y_3$ survive as generators in $H_{*}(\Omega
(F_4/T^4);\Z )$. If $y_2$ were not a generator in $H_{*}(\Omega
(F_4/T^4);\Z )$, we would have $y_2=\alpha x_1^4 +
\sum\limits_{i=2}^{4}\alpha _{i}x_1^3x_i+\sum\limits_{2\leq i<j\leq
4}\alpha _{ij}x_1^2x_ix_j+\beta x_1x_2x_3x_4$ for some $\alpha
,\alpha _{i},\alpha _{ij},  \beta \in \Z$. On the other hand, the
relation $2y_2=y_1^2$ from $H_{*}(\Omega F_4;\Z )$ becomes
$2y_2=x_1^4$ in $H_{*}(\Omega (F_4/T^4);\Z )$. This implies that
$2\alpha =1$ which is impossible. In the similar way we prove that
$y_3$ is also a generator in $H_{*}(\Omega (F_4/T^4);\Z )$. If it
were not, we would have $y_3=\alpha x_1^6 + \sum\limits_{i=2}^{4}
\alpha _{i}x_1^5x_i + \sum\limits _{2\leq i<j\leq 4}x_1^4x_ix_j+
\beta x_1^3x_2x_3x_4+\delta x_1^2y_2+\sum\limits _{1\leq i<j\leq
4}\delta _{ij}x_ix_jy_2$. From $H_{*}(\Omega F_4;\Z )$, we also have
that $3y_3=x_1^2y_2$. This together leads to $3\delta =1$ which is
impossible.
\end{proof}
\subsection{The integral homology of $\Omega (E_6/T^6)$}
\begin{remark}
The integral homology algebra $H_{*}(\Omega E_6;\Z )$ is described
in~\cite{M} and it is given by
\[
H_{*}(\Omega E_6;\Z )\cong \Z
[y_1,y_2,y_3,y_4,y_5,y_7,y_8,y_{11}]/\langle y_1^2-2y_2,
y_1y_2-3y_3\rangle \;,
\]
where $\deg y_{i}=2i$ for $i=1,2,3,4,5,7,8,11$.
\end{remark}
Using the same argument as for the previous cases, we deduce the
integral Pontrjagin homology of the based loop space on $E_6/T^6$.
\begin{theorem}\label{intE6/T}
The integral Pontrjagin homology ring of $\Omega (E_6/T^6)$ is given
by
\[
H_{*}(\Omega (E_6/T^6); \Z )\cong
\]
\[
\left(T(x_1,x_2,x_3,x_4,x_5,x_6)\otimes \Z
[y_1,y_2,y_3,y_4,y_5,y_7,y_8,y_{11}]\right)/I
\]
where $I= \langle x_k^2=x_px_q+x_qx_p=12y_1 \text{ for } 1\leq
k,p,q\leq 6,\; 2y_2=x_1^4,\; 3y_3=x_1^2y_2 \rangle$ and where $\deg
x_i=1$ for $1\leq i\leq 6$, and $\deg y_i=2i$ for
$i=2,3,4,5,7,8,11$.
\end{theorem}

%%% The bibliography %%%
\bibliographystyle{amsalpha}

\end{document}